# Can Laplace's formula model a deterministic universe that is *irreducibly* probabilistic?

Bhupinder Singh Anand

*(This paper includes, and extends, the concepts introduced in "Is a deterministic universe logically consistent with an irreducibly probabilistic Quantum Theory?")*

If we assume the Thesis that any classical Turing machine T, which halts on every $n$-ary sequence of natural numbers as input in a determinate time $t(n)$, determines a PA-provable formula, whose standard interpretation is an $n$-ary arithmetical relation $f(x_1, ..., x_n)$ that holds if, and only if, T halts, then we can define Laplace's formula recursively such that it can express the state of a deterministic quantum universe which is irreducibly probabilistic.

## Contents

# 1.  Introduction

Are our current theories of physics consistent with the concept of a universe that is completely deterministic, yet not pre-destined?

In other words, can the initial conditions and all physical laws at any instant, say, for instance, at the time of a projected Big Bang, be knowable completely in a manner that is consistent with our current theories of physics?

## 1.1  Is our universe pre-destined?

As Ian Stewart remarks ([St97], p6), the pre-quantum, classical concept was that our universe is both deterministic and pre-destined, a view expressed "eloquently" by Pierre Simon de Laplace, one of the leading mathematicians of the $18^{th}$ century, in his "*Philosophical Essays on Probabilities*":

> An intellect[1] which at any given moment knew all the forces that animate Nature and the mutual positions of the beings that comprise it, if this intellect were vast enough to submit its data to analysis, could condense into a single formula[2] the movement of the greatest bodies of the universe and that of the lightest atom: for such an intellect nothing could be uncertain; and the future just like the past would be present before its eyes.

---

[1] For convenience, we shall refer to such an intellect as Laplace's intellect "Li", pronounced "Lee".

[2] For convenience, we shall refer to such a formula as Laplace's formula "Lf", pronounced "Elf".

Such a notion of a pre-destined universe evolved as "the revolution in scientific thought that culminated in Newton led to a vision of the universe as some gigantic mechanism ... In such a vision, a machine is above all predictable" ([St97], p6).

## 1.2 Is our universe deterministic?

As Stewart observes further ([St97], p329), the post-quantum belief that our universe may be deterministic in a yet unknown, but fundamental, way (which may not necessarily be pre-destined) is reflected in Einstein's well-known remark, in the following excerpt from a letter to Max Born:

> You believe in the God who plays dice, and I in complete law and order in a world which objectively exists, and which I, in a wildly speculative way, am trying to capture. I firmly *believe*, but I hope that someone will discover a more realistic way, or rather a more tangible basis than it has been my lot to do. Even the great initial success of the quantum theory does not make me believe in the fundamental dice game ... .

## 1.3 Is quantum mechanics 'irreducibly probabilistic'?

Stewart notes, further ([St97], p330), the prevalent view that, despite Einstein's predilections, the universe, or at least our present quantum mechanical description of it, is of an "irreducibly probabilistic character". He then suggests that we may need to seriously consider the "... possibility of changing the theoretical framework of physics altogether, replacing quantum uncertainty by deterministic chaos, as Einstein would have liked".

## 1.4  Can we mathematically express a deterministic, yet quantum, universe that is not pre-destined?

Clearly, the extent to which we can address this issue meaningfully will depend on how precisely we can express the above concept, of a deterministic, yet quantum, universe that is not pre-destined, within a constructive[3], and intuitionistically unobjectionable language.

In this paper, we broadly address the above issue of whether such expression can be made mathematically for a Big Bang type of universe U that, we assume, creates at most denumerable particles at the instant of creation, all super-posed upon each other at a single "point", where the states, and inter-actions, of the particles are dynamic and follow deterministic laws, and where the future can reasonably be interpreted as *irreducibly* probabilistic.

Specifically, we consider the following:

(*i*)  Can we constructively define intuitive concepts such as "particle", "deterministic", "quantum", "knowable", "pre-destined", and, most significantly, "Laplace's formula" mathematically?

(*ii*)  Is there a mathematical description of the state of a given particle $P$ at a given inter-action[4] $t$, such that, given any property $x$ of P, there is always an effective method by which we can determine the value $y$ of the property $x$ uniquely?

---

[3] We term a mathematical language constructive if, and only if, there is an effective method of assigning truth-values to the propositions of the language (cf. [An02c], §5).

[4] The term "inter-action", as used here, and henceforth, should be taken to correspond intuitively to the term "instant". However, such instants are only those discrete moments of time when the particle "inter-acts" with some other particle in the universe. This presumes that the particle's inter-action with the universe is not continuous, and so the state, or signature, of the particle is known only as, and changes only as, a result of an inter-action.

(*iii*) Can such a description allow for the possibility that, given any arbitrary set of properties of a given particle *P*, there is no effective method such that we can always determine the precise values of the given properties simultaneously at any inter-action *t*?

(*iv*) Can such a description allow for the possibility that the sequence of values of some property *x* of a given particle *P*, measured at discrete inter-actions $t_1$, $t_2$, ..., may be Li-unpredictable?

**Preliminary Definition 1**: A class of natural numbers is Li-predictable if, and only if, it is a well-defined mathematical object[5] (such as, for instance, any well-defined set determined by a consistent set theory[6]). We assume that, in this case, Li is aware of any member of the class directly, without depending on any effective method for computing it.

**Preliminary Definition 2**: A sequence of natural numbers is Li-unpredictable if, and only if, it does not define a mathematical object. We assume that, in this case, Li is not aware of the value of the *i*'th term directly for any given *i*, but must depend on an effective method for computing it.

(*v*) Can we assume that fundamental particles have a distinct identity throughout their "lifetime" in such a description?

---

[5] The formal definition of a mathematical object is given in Anand ([An02c], §1.2, Definition 1(*iii*)).

[6] We note that, since there are recursive number-theoretic relations and functions that cannot be introduced through definition into any formal system of Arithmetic without inviting inconsistency ([An02c], Corollary 1.1), such a theory cannot contain a standard Axiom Schema of Separation (Comprehension) such as:

Suppose X is a set and f is a formula then $\{x$ is in $X$: f $(x)\}$ is also a set. More formally, (A$x$)(E$y$)($z$ is in $y <=> (z$ is in $x$ and f $(x)$).

(*vi*) Can two distinctly different particles "occupy" the same space simultaneously in such a description?

## 1.5  Quantum concepts in a constructive framework of recursive functions and Peano's Arithmetic

The underlying thesis of this paper is that a constructive, non-classical, mathematical language for quantum mechanical descriptions, which addresses these questions appropriately, can, indeed, be built around asymmetric ([An02c], §2.1) recursive number-theoretic relations. These relations, characteristically, are not the standard interpretations of any of their formal representations in any standard first order Arithmetic that is based on Dedekind's formulation of the Peano Axioms for the natural numbers [An02c].

## 1.6  A possible ambiguity in classical Turing computability

However, in order to address these questions appropriately, we argue that there is a need to make the concept of classical Turing-computability[7] less ambiguous. This ambiguity is highlighted when we consider the two meta-theses:

(*a*)  **Classical Halting Thesis**: Every PA-formula[8], that is a true[9] $n$-ary arithmetical[10] relation $f(x_1, \ ..., x_n)$ under the standard interpretation[11], determines a classical

---

[7] We follow Mendelson's definition of a classical Turing machine ([Me64], p229), and of classical Turing-computability ([Me64], p231).

[8] By "PA", we mean a classical, standard, first order formal system of Arithmetic such as Mendelson's formal system S ([Me64], p102). We note, however, that the arguments of this paper would also hold for Gödel's formal system P ([Go31a], p9).

In the definitions of logical terms, we generally follow Mendelson's definitions [Me64] such as, for instance:

(*i*)  A formula of a formal system S is any finite sequence of the alphabet (assumed to consist of a countable set of symbols) of S ([Me64], p29).

Turing machine T that halts on every $n$-ary sequence $k_1$, ..., $k_n$ of natural numbers as input if, and only if, $f(k_1, ..., k_n)$ holds.

(*b*) **Quantum Halting Thesis**: Any classical Turing machine T that halts on every $n$-ary sequence of natural numbers as input in a determinate time $t(n)$ determines a PA-provable formula, whose standard interpretation is an $n$-ary arithmetical relation $f(x_1, ..., x_n)$ such that, for any sequence $k_1$, ..., $k_n$ of natural numbers, $f(k_1, ..., k_n)$ holds if, and only if, T halts.

Clearly, the two meta-theses are mutually inconsistent, since (*a*) would determine a classical Turing machine T that halts on every $n$-ary sequence $k_1$, ..., $k_n$ of natural numbers as input even for a PA-unprovable formula that is a true $n$-ary arithmetical relation $f(x_1, ..., x_n)$ under the standard interpretation[12].

---

(*ii*) A well-formed formulas of S is any one of a specially selected sub-set of the formulas of S, where the selection is such that, given any formula $F$, there is always an effective method to decide whether $F$ is a well-formed formula of S or not ([Me64], p29).

(*iii*) A formula of S is S-provable if, and only if, there is a finite proof sequence of S-formulas each of which is either an axiom of S, or an immediate consequence of two of the preceding formulas in the sequence by the rules of inference of S ([Me64], p29).

[9] We follow Tarski's definitions of the "satisfiability" and "truth" of a formula under a given interpretation ([Me64], p49-53).

We note that Tarski's definitions implicitly assume that a formula [$F$] is true under an interpretation if, and only if, there is an effective method (native to the interpretation) for determining that the formula [$F$] is satisfied under the interpretation for all sequences in the domain of the interpretation ([An02c], §II(5.1)). It follows by Church's Thesis ([Me64], p147, footnote) that, if the PA-formula [$F$] is true under the standard interpretation, then it is recursive.

[10] We follow Mendelson's definition of an arithmetical relation (predicate) ([Me64], p135, Ex. 2a). We note that, in standard PA, this definition is equivalent to Gödel's definition of arithmetical relations ([Go31a], p29).

[11] We follow Mendelson's definitions of interpretation ([Me64], p49) and standard interpretation ([Me64], p107).

[12] In his classic 1931 paper [Go31a], Gödel constructs such a formula in an intuitionistically unobjectionable manner.

Further, since the conclusion in (*a*), and the premise in (*b*), are not effectively verifiable, the two meta-theses are essentially undecidable. Prima facie, we thus have two, apparently equally sound, but essentially different, systems of Arithmetic based on constructive, and intuitionistically unobjectionable[13], interpretations of standard PA; reasonably, they should have significantly different consequences.

**1.7 Three theorems**

We argue that this is, indeed, so, in the following:

**Theorem 1**: The principle of Quantum Uncertainty is inconsistent with the Classical Halting Thesis.

**Theorem 2**: The principle of Quantum Uncertainty is logically consistent with the Quantum Halting Thesis.

**Theorem 3**: The Quantum Halting Thesis implies that every partial recursive function can be effectively defined as a total[14] function.

## 2. Defining a deterministic universe

We start by introducing the critical concept of a formally undecidable proposition, first constructed by Gödel in his seminal 1931 paper, as the basis for defining a mathematical function that, loosely speaking, corresponds to the classical quantum state of a particle.

---

[13] We use the terms "constructive" and "intuitionistically unobjectionable" in the sense in which Gödel used them to describe his formal reasoning in his seminal 1931 paper [Go31a].

[14] Loosely speaking, a number-theoretic function is total if, and only if, it is defined for every set of natural number values of its free variables.

**Definition 1**: The standard interpretation of a unary PA-formula is a quantum signature if, and only if, the formula is PA-unprovable but true under the standard interpretation.

The significance of a quantum signature $[F(x)]$[15] is that, given any natural number $k$, we can always construct an individually[16] effective[17] method for determining that $F(k)$ holds, but there may be no uniformly effective method, even for Laplace's intellect, such that, given any $k$, we can always determine that $F(k)$ holds.

Treating the terms "particle", "property" and "value" as undefined primitives, we intuitively assert the:

**Expressibility Thesis**: Any physical property of any given particle that can be measured by some effective method can be expressed by a finite string in a suitably defined, recursively enumerable, language[18] L, which has a finite alphabet. Similarly,

---

[15] We denote the PA-formula whose interpretation is $F$ by $[F]$.

[16] The thesis underlying this paper is based on the argument that we can differentiate between intuitive concepts such as "individually decidable", "individually computable", "individually terminating routine" etc., and their corresponding, and contrasting, intuitive concepts such as "uniformly decidable", "uniformly computable", "uniformly terminating routine" etc. ([An02c], §II(5)).

Loosely speaking, we would use the terminology, for instance "individually decidable", if, and only if, say, a PA-formula with free variables is decidable individually, as satisfiable or not, in an interpretation M for any given instantiation. This corresponds to, but is not necessarily a consequence of, the case where every instantiation of the formula is PA-provable. We note that the closure of the formula may not be P-provable

In contrast, we would use the terminology "uniformly decidable" if, and only if, a PA-formula with free variables is decidable jointly, as true or not, in the interpretation M as an infinitely compound assertion of all of its instantiations. This corresponds to, and is a consequence of, the case where the closure of the formula is PA-provable.

The concept of an intuitively terminating routine may be taken as corresponding to the classical concept of an "effective procedure" ([Me64], p207).

[17] We broadly follow Mendelson's approach to "effectiveness" ([Me64], p207).

[18] We define a language L as recursively enumerable if there exists a Turing machine that accepts every string of the language, and does not accept strings that are not in the language.

any value of such a physical property that can be determined by some effective method can also be expressed by a finite string in L.

It immediately follows that:

**Lemma 1**: Any property can be uniquely Gödel-numbered[19] by a natural number, and any value of a property can, similarly, be uniquely Gödel-numbered by a natural number.

We note that not every natural number need correspond to a property, and that not every natural number need correspond to the value of a property. We also note that, even if the number of all possible particles, properties and values in a deterministic universe U are not denumerable, since we can only express denumerable strings within L, in this paper we only consider the recursively enumerable set of properties for which we have an effective method of measurement.

We next define:

**Definition 2**: The signature of a given particle $P_i$ at any inter-action $t$ is a number-theoretic relation $p(i, t, x, y)$ that holds if, and only if, for any given property represented by $x$, there is a unique natural number $y$, where the domain of $x$ is assumed to be a recursively enumerable set[20] of the Gödel-numbers of all the properties that are expressible in L[21], and the domain of $y$ is assumed to be a recursively enumerable set of the Gödel-numbers of all values that are expressible in L.

---

[19] By "Gödel-numbered", we mean the assignment of a unique natural number to the formal expression of any concept within a language ([Me64], p135).

[20] " 'A set of natural numbers is called recursively enumerable (r.e) if, and only if, it is either empty or is the range of a recursive function. Intuitively, if we accept Church's Thesis, then a recursively enumerable set is a collection of natural numbers which is generated by some mechanical process." ([Me64], p250)

[21] Prima facie, it is not obvious whether the properties need to be mutually independent, where we define a set of properties as mutually independent if, and only if, any property in the set can assume any

We can now define:

**Definition 3**: The history $H_{i,\,t}$ of any given particle $P_i$[22] at any specified inter-action $t$[23] is the set of particle signatures $\{p(i, s, x, y): 0 =< s =< t\text{-}1\}$.

**Definition 4:** For any given particle $P_i$, we define a given property[24] $k$ as deterministic if, at any given inter-action $t$, there is an effective method[25] to determine a unique value $m_{k,\,t}$ such that $p(i, t, k, m_{k,\,t})$ holds.

**Definition 5**: A given universe U is deterministic if, for any given particle $P_i$, any given property $k$, and any given inter-action $t$, there is an effective method to determine a unique value $m_{k,\,t}$ such that $p(i, t, k, m_{k,\,t})$ holds.

We now note that:

---

permitted values in its range without affecting the values of any other property in the set. Formally, if the properties are denoted by $k_1, k_2, ...,$ and the permitted range of values of the property denoted by $k_i$ is $R_i$, then the set of properties is mutually independent if, and only if, every permissible sequence of values, denoted by $m_1, m_2, ...,$ where $m_i$ is in $R_i$, is a valid set of values for the set of properties denoted by $k_1, k_2, ....$

[22] Although we cannot assume that a given history necessarily defines a unique particle, for the moment, it is convenient, loosely speaking, to view the particle $P$ as a Platonic mathematical entity with a unique identity. An entity, moreover, that manifests itself only through its history - in other words, through its set of signatures at its various inter-actions during its lifetime.

[23] We note that $t$ here is a natural number that is assumed to be the $t$'th inter-action of the particle $P_i$, with some other particle in the universe, with reference to an initial signature $p(i, 0, x, y)$. In other words, we assume that, for any property $k$, the values $p(i, 0, k, m_{k\,0}), p(i, 1, k, m_{k,\,1}), ..., p(i, t\text{-}1, k, m_{k,\,t\text{-}1})$ are already known.

[24] Where there is no obvious ambiguity, we shall henceforth use the terms "property $k$" and "value $m_{k,\,t}$" when we actually mean "property whose Gödel-number is $k$", and "value whose Gödel-number is $m_{k,\,t}$ at inter-action $t$".

[25] We note that the concept of an "effective method" is essentially intuitive, and, loosely speaking, may be taken as corresponding to the result of some experiment that is based on assuming deterministic, fundamental, Laws of Nature in the classical, Newtonian sense as, apparently, intended by Einstein in his above cited letter.

**Lemma 2**: If the properties expressible in L are recursively enumerable, and we assume Turing's Thesis[26], then we can effectively determine whether $p(i, t, k, m)$ holds or not in a given deterministic universe U for any given natural numbers $i, t, k, m$.

*Proof*: By definition, and the Turing Thesis, there is a classical Turing machine T such that, given any natural numbers $i, t, k, m$ as input, T halts if, and only if, $p(i, t, k, m)$ holds. If $k$ is not the Gödel-number of a property, then $p(i, t, k, m)$ is undefined, and we assume that, in such a case, T will loop on the input $k$.

We note that the definition of any classical Turing machine T can be extended to allow for recording, in its infinite memory[27], of every instantaneous tape description[28] during the operation of a classical Turing machine.

Now, T can always be meta-programmed to compare its current instantaneous tape description with the finite set of previous instantaneous tape descriptions during its current operation, and to self-terminate if an instantaneous tape description repeats itself, i.e. at the onset of a loop[29]. We can, thus, use this condition to effectively define $p(i, t, k, m)$ as holding vacuously if $k$ is not the Gödel-number of a property.

---

[26] We take Church's Thesis to essentially state that a number-theoretic function is effectively computable if, and only if, it is recursive, and a partial number-theoretic function is effectively computable if, and only if, it is partial recursive ([Me64], p227). We note too, that a number-theoretic function is Turing-computable if, and only if, it is partial recursive ([Me64], p233, Corollary 5.13 & p237, Corollary 5.15). We thus have the Turing Thesis [Tu36] that a number-theoretic function is effectively computable if, and only if, it is Turing-computable.

[27] We may visualise T as a virtual teleprinter that copies, and records, every instantaneous tape description in an area of the machine that acts as a dynamic, random access, memory during any operation of T.

[28] Cf. Mendelson ([Me64], p230).

[29] Cf. Anand ([An02c], §II(7.1)).

Hence we can effectively determine whether $p(i, t, k, m)$ holds or not for any given natural numbers $i, t, k, m$.¶[30]

We thus have:

**Lemma 3**: If we assume the Church Thesis, then the number-theoretic relation $p(i, t, x, y)$, describing the signature of the particle $P_i$ at inter-action $t$, is recursive[31] in a given deterministic universe U.

*Proof*: This follows from the argument that every Turing-computable function $f$ is recursive if it is total[32].

We note that, in view of Lemma 3, the particle signature $p(i, t, x, y)$ is essentially Laplace's formula for a given deterministic universe U, since it essentially condenses "into a single formula the movement of the greatest bodies of the universe and that of the lightest atom", and "all the forces that animate Nature and the mutual positions of the beings that comprise it".

We also note that the particle signature $p(i, t, x, y)$, at any inter-action $t$ in a given deterministic universe U, can be viewed as a point in an infinite dimensional Hilbert space, whose co-ordinates are $(m_{t, 1}, m_{t, 2}, ...)$, where $m_{t, k}$ is the value along the $k$ axis that satisfies $p(i, t, k, m_{t, k})$. We can thus define:

**Definition 6**: A Gödelian point $g_i(t)$ of the particle $P_i$ is a point in an infinite dimensional Hilbert space, whose co-ordinates are $(m_{t, 1}, m_{t, 2}, ...)$, where $m_{t, k}$ is the value along the $k$ axis that satisfies $p(i, t, k, m_{t, k})$.

---

[30] We can thus treat the particle signature $p_t(x, y)$ at any inter-action $t$ as a point in an infinite dimensional Hilbert space, whose co-ordinates are $m_1, m_2, ...$, where $m_k$ is the value of the property $k$..

[31] We follow Mendelson's ([Me64], p120) definitions of recursive functions and relations.

[32] Cf. Mendelson ([Me64], p233, Corollary 5.13).

**Definition 7**: The Gödelian world-line[33] G($i$) of the particle $P_i$ is the non-terminating sequence of particle signatures $p(i, 0, x, y)$, $p(i, 1, x, y)$ ... in a given deterministic universe U.

**Definition 8**: A Gödelian world-line is a fundamental particle $P_i$ in a given deterministic universe U if, and only if, the non-terminating sequence of the Gödelian points, $g_i(t)$, t >= 0, is a well-defined mathematical entity.

## 2.1  Defining a quantum universe

We next define:

**Definition 9**: The PA-formula [$Q_t(x, y)$], which expresses[34] the recursive relation $p_i(t, x, y)$ in PA, is a quantum signature if, and only if, the formula [(E!$y$)$Q_t(x, y)$][35] is PA-unprovable, but true under the standard interpretation of PA.

We note that, by definition, if the PA-formula [$Q_t(x, y)$] expresses the recursive relation $p_t(x, y)$, then, for any given any natural numbers $k, m$:

(*i*)  (E!$y$)$Q_t(k, y)$ holds, and

(*ii*)  $Q_t(k, m)$ holds if, and only if, $p_t(k, m)$ holds.

We can now define:

---

[33] A Gödelian world-line is thus a step function in a Hilbert-space ([Pe90], p332; [Ru53], p229), p219). Treating this as defining a unique particle would correspond, loosely speaking, to the classical relativistic representation of a particle as a continuous world-line in Minkowski space-time ([Pe90], p250).

[34] We follow Mendelson's ([Me64], p117-8) definitions of expressibility and representability of number-theoretic functions.

[35] The notation "(E!$x$)" denotes uniqueness.

**Definition 10**: A Gödelian world-line is quantum if the PA-representation of its signature at every inter-action is a quantum signature.

**Definition 11**: A deterministic universe is pre-destined if it does not contain any quantum Gödelian world-line.

**Definition 12**: A deterministic universe is partially pre-destined if, and only if, it is not pre-destined.

**Definition 13**: A partially pre-destined universe is quantum if, and only if, it contains only quantum Gödelian world-lines.

**Lemma 4**: A quantum universe Q contains no fundamental particles.

## 3.  The Classical Halting Thesis in a quantum universe

We now have that:

**Lemma 5**: If we assume the Classical Halting Thesis (*a*) in a deterministic universe, then, for any given particle $P$, there is a classical Turing machine T such that, for any input $k$, T will halt and return the value $m$ if, and only if, $p_t(k, m)$ holds.

*Proof*: T will halt on input $k$ if, and only if, $(E!y)Q_t(k, y)$ holds. By our hypothesis, for any given $k$, $(E!y)Q_t(k, y)$ holds uniquely for some value $m$. Hence, for given input $k$, T must run a sub-routine that will halt on, and can be programmed to return, $m$.¶

This would clearly imply that:

**Lemma 6**: If we assume the Classical Halting Thesis (*a*), then there is some effective method that can completely determine the values of any finite number of, arbitrarily selected, properties simultaneously for any given particle in a deterministic universe.

We define:

**Definition 14**: The principle of Quantum Uncertainty holds in a universe if, and only if, there are properties of some particles whose values cannot be completely determined simultaneously.

It follows that:

**Theorem 1**: The principle of Quantum Uncertainty is inconsistent with the Classical Halting Thesis.

## 4.   The Quantum Halting Thesis in a quantum universe

However, in a quantum universe, the formula $[(E!y)Q_t(x, y)]$ is not PA-provable for any given particle. It follows that:

**Lemma 7**: If we assume the Quantum Halting Thesis (*b*) in a deterministic universe, then, for any given particle $P$, there is no classical Turing machine T such that, for any given input $k$, T will halt and return the value $m$ if, and only if, $p_t(k, m)$ holds.

This, now, clearly implies that:

**Lemma 8**: There is no effective method that can always determine the values of a finite number of, arbitrarily selected, properties of any given particle $P$ simultaneously, even though, by our premise, given any property $k$, there is always some effective method that will completely determine its value $m_k$ at inter-action $t$.

We thus conclude that:

    **Theorem 2**: The principle of Quantum Uncertainty is consistent with the Quantum Halting Thesis[36].

## 5. Some consequences

We note that $t$ in *Definition 2* is a natural number that is assumed to be the $t$'th inter-action of the particle $P$ with some other particle in the universe, necessarily with reference to some initial signature $p_0(k, m_{k, 0})$[37]. In other words, we assume that, for any given particle P and any given property $k$, the signature $p_0(k, m_{k, 0})$, $p_1(k, m_{k, 1})$, ..., $p_{t-1}(k, m_{k, t-1})$, ... is assumed known to Laplace's deity.

### 5.1 The Gödel ß-function

We note further that, for any given property $k$ and inter-action $t = t_n$, we can define a recursive Gödel $\beta$-function[38] $q_{P, k}(b, c, i)$, such that $q_{P, k}(b, c, i) = m_{k, i}$ for all $0 =< i =<$

---

$n$[39], where $m_{k,\,i}$ is the value of the property $k$ at the inter-action $t_i$, and $b$, $c$ are natural numbers that can be constructively determined by the sequence $m_{k,\,0}$, $m_{k,\,1}$, ..., $m_{k,\,n}$.

## 5.2 Properties with random values in a deterministic universe

If we assume that the inter-actions $t_0$, $t_1$, ..., $t_n$ refer to the series of inter-actions that occur during a particles lifetime[40], then we may reasonably define:

**Definition 15**: A property $k$ is predictable if the sequence of values represented by the sequence $m_{k,\,0}$, $m_{k,\,1}$, ..., $m_{k,\,n}$, as the number $n$ of inter-actions increases indefinitely, is a Cauchy sequence; otherwise $k$ is random.

Prima facie, it thus follows that:

**Lemma 9**: Both predictable and random properties are consistent with a deterministic universe.

## 5.3 Probability of an inter-action yielding a given value

We note, further, that, for any given sequence $m_{k,\,0}$, $m_{k,\,1}$, ..., $m_{k,\,n}$, there are denumerable pairs of natural numbers $(b, c)$ that define Gödel $\beta$-functions relative to the sequence; however, although the functions will all yield identical values of $m_{k,\,i}$ for $0 =< i =< n$, each of them will yield different sequences $m_{k,\,n+1}$, $m_{k,\,n+2}$, ....

The challenge of any theory, thus, would be to determine, firstly, which of these functions are such that all the terms of the non-terminating sequence are Gödel-numbers of values of the property $k$; secondly, the probability that the result $m_{k,\,n+1}$ of the inter-action $t_{n+1}$, for any pair $(b, c)$, will represent any given value of the property $k$.

---

[39] A notation such as "=<" is to be read naturally as "equal to or less than".

[40] Since, say, the Big Bang.

We thus have that:

**Lemma 10**: A sufficient condition for a deterministic universe to be pre-destined is that every property of every particle P have a characteristic Gödel $\beta$-function that determines the value of the property at any inter-action $t_n$ for all $n >= 0$.

**Lemma 11**: A necessary condition for a deterministic universe to be pre-destined is that the initial signature $p_0(k_0, m_0)$, $p_0(k_1, m_1)$, ..., $p_0(k_{t-1}, m_{t-1})$, ... of every particle P is assumed known.

If we now define:

**Definition 16**: A universe is predictable if there is always an effective method such that, for any given particle $P$, any given deterministic property $k$, and any given inter-action $t$, we can determine a unique value $m_k$ such that $p_t(k, m_k)$ holds.

We then have that:

**Lemma 12**: A sufficient condition for a deterministic universe to be unpredictable is that the initial signature $p_0(k_0, m_0)$, $p_0(k_1, m_1)$, ..., $p_0(k_{t-1}, m_{t-1})$, ... of every particle P is assumed unknowable.

## 5.4 Superposition

We also note that, if $q_{P, k}(b, c, i) = m_{k, i}$ for all $0 =< i =< n$, where $m_{k, i}$ is the value of the property $k$ at the inter-action $t_i$, and $q'_{P, k}(b', c', i) = m'_{k, i}$ for all $0 =< i =< n'$, where $m'_{k, i}$ is the value of the property $k$ at the inter-action $t'_i$, we can always combine the two sequences of inter-actions $m_{k, 0}$, $m_{k, 1}$, ..., $m_{k, n}$, and $m'_{k, 0}$, $m'_{k, 1}$, ..., $m'_{k, n'}$, appropriately[41] to yield a sequence $m''_{k, 0}$, $m''_{k, 1}$, ..., $m''_{k, n''}$, where $n'' = n + n'$.

---

[41] We assume that there is an effective way of ensuring that the combined sequence $t''_i$, $0 =< i =< n''$, is also chronological.

We thus have a Gödel $\beta$-function $q''_{P,k}(b'', c'', i) = m''_{k,i}$ for all $0 =< i =< n''$, where $m''_i$ is the value of the property $k$ at the inter-action $t''_i$, which may be considered as the superposition of the Gödel $\beta$-function $q'_{P,k}(b', c', i)$ on $q_{P,k}(b, c, i)$.

It follows that:

**Lemma 13**: The probability[42] that $q''_{P,k}(b'', c'', n''+1) = m$ for a given $m$ is always equal to, or greater than, the probability that $q_{P,k}(b, c, n+1) = m$.

*Proof*: The sequence $m_0, m_1, ..., m_n$ is a proper sub-sequence of $m''_0, m''_1, ..., m''_{n''}$. Hence, if $q''_{P,k}(b'', c'', n''+1) = m$, then there is always some pair $(b, c)$ such that $q_{P,k}(b, c, n''+1) = m$. The converse is not true.

Hence the set of Gödel $\beta$-functions for which $q''_{P,k}(b'', c'', i) = m''_i$ for all $0 =< i =< n''$, where $m''_i$ is the value of the property $k$ at the inter-action $t''_i$, is a proper sub-set of the set of Gödel $\beta$-functions for which $q_{P,k}(b, c, i) = m_i$ for all $0 =< i =< n$, where $m_i$ is the value of the property $k$ at the inter-action $t_i$.¶

## Appendix 1: Defining every partial recursive function as total

The above consequences emerged incidentally out of the following, prima facie unrelated, argument, which effectively defines every partial recursive function as total.

**Theorem 3**: The Quantum Halting Thesis implies that every partial recursive[43] number-theoretic function $f(x_1, ..., x_n)$ is constructively definable as a total function.

---

*Proof*: We assume that *f* is obtained from the recursive function *g* by means of the unrestricted *μ*-operator, so that $f(x_1, ..., x_n) = \mu y(g(x_1, ..., x_n, y) = 0)$.

Given any sequence $k_1, ..., k_n$ of natural numbers, we can then define:

(*i*)   a classical Turing machine T1 that will halt if, and only if, $(g(k_1, ..., k_n, k) = 0)$ holds for some natural number k;

(*ii*)   a classical Turing machine T2 that will halt if, and only if, ~$(g(k_1, ..., k_n, k) = 0)$ holds for some natural number k;

(*iii*)   a classical Turing machine T3[44] that will halt if, and only if, $[H(k_1, ..., k_n, y)]$ is PA-provable, where $[H(x_1, ..., x_n, y)]$ is the PA-formula that represents $(g(x_1, ..., x_n, y) = 0)$ in PA.

If we now assume the Quantum Halting Thesis, it follows that:

(*a*)   if T1 halts, without looping, on all natural numbers, then $[H(k_1, ..., k_n, y)]$ is PA-provable;

(*b*)   similarly, if T2 halts, without looping, on all natural numbers, then $[H(k_1, ..., k_n, y)]$ is PA-provable.

Hence, if $(g(k_1, ..., k_n, y) = 0)$ holds for all natural numbers y, then:

(*c*)   either $[H(k_1, ..., k_n, y)]$ is PA-provable, and T3 will halt;

---

the least number *k* (if such exists) such that, if $0=<i=<k$, $g(x_1, ..., x_n, i)$ exists and is not 0, and $g(x_1, ..., x_n, k) = 0$. We note that, classically, *f* may not be defined for certain *n*-tuples; in particular, for those *n*-tuples $(x_1, ..., x_n)$ for which there is no *y* such that $g(x_1, ..., x_n, y) = 0$.

[44] The definition of T3 is based on Gödel's recursive number-theoretic relation *xBy* ([Go31a], Definition 45), which holds if, and only if, *x* is the Gödel-number of a PA-proof of the PA-formula whose Gödel-number is *y*.

(*d*)   or, if [$H(k_1, ..., k_n, y)$] is not PA-provable, then T2 will loop, and self-terminate, for some natural number $k$.

If we now run T1, T2 and T3 simultaneously, then we are assured that at least one of them will halt, or self-terminate, for some natural number $k$. We can thus use the condition indicated by the halting, or self-terminating, state to effectively define the partial recursive function $f(x_1, ..., x_n)$ appropriately as a total partial recursive function.¶

## Appendix 2: Collapse of the wave function

Consider this:

(*i*) Gödel has proved that there is an arithmetic formula [$G(x)$] such that, for any given $k$, [$G(k)$] is provable.

(*ii*) Hence, for any given $k$, there is always some effective method for evaluating the arithmetic expression $G(k)$.

(*iii*) Gödel has also proved that [(A$x$)$G(x)$] is not provable.

(*iv*) *Thesis*: There is no uniform effective method (algorithm/Turing machine) that can evaluate the arithmetic expression $G(x)$ for any given $x$.

(*v*) Thus, $G(n)$ is individually computable, but not uniformly computable.

(*vi*) *Theorem (provable by induction)*: For any given $k$, we can always find some effective method (algorithm/Turing machine) $T(k)$ that can compute $G(n)$ for all $n<k$, i.e. $T(k)$ terminates for all $n<k$, but it "loops" on input $k$. (Note: All methods that evaluate $G(n)$ for all $n<k$ cannot be non-terminating on input $k$; this would imply that $G(k)$ is undefined, which contradicts (*ii*).)

(*vii*) *Quantum interpretation*: The process of finding $T(k+1)$ can be corresponded, firstly, to the act of finding a suitable method of measuring the value $G(k)$ precisely, and, secondly, to the collapse of the wave function at $k$ as a result of the measurement; we then have the new "state" $T(k')$, which can evaluate the value of $G(n)$ for all $n<k'$, where $k<k'$, but not beyond!

*Author's e-mail: anandb@vsnl.com*

        (*Updated: Monday $16^{th}$ June 2003 3:21:21 PM by re@alixcomsi.com*)